\documentclass{fundam}

\usepackage{url} 
\usepackage[ruled,lined]{algorithm2e}
\usepackage{graphicx}
\usepackage{tikz}
\allowdisplaybreaks
\usetikzlibrary{automata, positioning, arrows}

\begin{document}


\setcounter{page}{61}
\publyear{22}
\papernumber{2130}
\volume{187}
\issue{1}

    \finalVersionForARXIV


\title{A Note of Generalization of Fractional ID-factor-critical Graphs}

\author{Sizhong Zhou\thanks{Address of correspondence: School of Science, Jiangsu University of Science and
                        Technology, Changhui Road 666, Zhenjiang, Jiangsu 212100, China.  \newline \newline
          \vspace*{-6mm}{\scriptsize{Received July 2019; \ accepted August  2022.}}}
          \\
School of Science\\
Jiangsu University of Science and Technology\\
Zhenjiang, Jiangsu 212100, China\\
zsz\_cumt{@}163.com
}

\maketitle

\runninghead{S. Zhou}{A Note of Generalization of Fractional ID-factor-critical Graphs}

\begin{abstract}
  In communication networks, the binding numbers of graphs (or networks) are often used to measure the vulnerability and
robustness of graphs (or networks). Furthermore, the fractional factors of graphs and the fractional ID-$[a,b]$-factor-critical covered
graphs have a great deal of important applications in the data transmission networks. In this paper, we investigate the relationship
between the binding numbers of graphs and the fractional ID-$[a,b]$-factor-critical covered graphs, and derive a binding number condition
for a graph to be fractional ID-$[a,b]$-factor-critical covered, which is an extension of Zhou's previous result [S. Zhou, Binding numbers
for fractional ID-$k$-factor-critical graphs, Acta Mathematica Sinica, English Series 30(1)(2014)181--186].
\end{abstract}

\begin{keywords}
network, graph, binding number, fractional $[a,b]$-factor, fractional ID-$[a,b]$-factor-critical covered graph.
\end{keywords}

\section{Introduction}

We investigate the fractional factor problem of graphs, which can be regard as a relaxation of the well-known cardinality matching
problem. It has wide-ranging applications in many distinct fields such as scheduling, network design, circuit layout, combinatorial
design and combinatorial polyhedron. For example, if we consider some large data packets to be sent to several distinct destinations
through some channels in a communication network, and to improve the efficiency of the network, then we may partition the large
data packets into small parcels. The feasible assignment of data packets can be considered as a fractional flow problem which is
also described as a problem of fractional factor in a graph.

In the process of data transmission, if some special nodes (i.e., nonadjacent nodes) are damaged and we require that a channel is
assigned, the possibility of data transmission in a communication network is considered as the existence of fractional
ID-factor-critical covered graph. Naturally, the existence of fractional ID-factor-critical covered graphs plays an important role
in data transmission networks. Several maturing methods on graph based network design were derived by de Araujo, Martins and Bastos
\cite{AMB}, Ashwin and Postlethwaite \cite{AP}, Fardad, Lin and Jovanovic \cite{FLJ}, Lanzeni, Messina and Archetti \cite{LMA},
Pishvaee and Rabbani \cite{PR}, and Rahimi and Haghighi \cite{RH}.

The graphs studied in this paper are simple. We denote a graph with vertex set $V(G)$ and edge set $E(G)$ by $G=(V(G),E(G))$.
For $x\in V(G)$, the set of vertices adjacent to $x$ in $G$ is said to be the neighborhood of $x$, denoted by $N_G(x)$, and
$|N_G(x)|$ is said to be the degree of $x$ in $G$, denoted by $d_G(x)$. Set $\delta(G)=\min\{d_G(x):x\in V(G)\}$. For any
$S\subseteq V(G)$, we write $N_G(S)=\bigcup\limits_{x\in S}N_G(x)$. The subgraph of $G$ induced by $S$ is denote by $G[S]$, and
$G-S=G[V(G)\setminus S]$. A vertex set $S\subseteq V(G)$ is called independent if $G[S]$ does not admit edges. Let $S$ and $T$ be
two disjoint subsets of $V(G)$. We denote by $e_G(S,T)$ the number of edges with one end in $S$ and the other end in $T$. The
binding number of $G$ is defined by
\begin{align*}
bind(G)=min\left\{\frac{|N_G(X)|}{|X|}:\emptyset\neq X\subseteq V(G),N_G(X)\neq V(G)\right\}.
\end{align*}

For two positive integers $a$ and $b$ with $a\leq b$, an $[a,b]$-factor of $G$ is a spanning subgraph $F$ of $G$ such that
$a\leq d_F(x)\leq b$ holds for all $x\in V(G)$. Let $h: E(G)\rightarrow [0,1]$ be a real-valued function from the edge set $E(G)$
to the real number interval $[0,1]$. If $a\leq\sum\limits_{e\ni
x}h(e)\leq b$ holds for any $x\in V(G)$, then we call $G[F_h]$ a fractional $[a,b]$-factor of $G$ with indicator function $h$,
where $F_h=\{e: e\in E(G), h(e)>0\}$. A fractional $[k,k]$-factor is simply called a fractional $k$-factor. A graph $G$ is fractional
ID-$[a,b]$-factor-critical if $G-I$ admits a fractional $[a,b]$-factor for any independent set $I$ of $G$. A fractional
ID-$[k,k]$-factor-critical graph is simply called a fractional ID-$k$-factor-critical graph. A great deal of results on the
topic with factors in graphs, fractional factors in graphs and fractional ID-factor-critical graphs can refer to Wang and Zhang
\cite{WZi,WZo,WZr}, Zhou, Sun and Bian \cite{ZSB}, Zhou \cite{Za,Za0,Zr}, Zhou and Bian \cite{ZB}, Haghparast and Kiani \cite{HK},
Hasanvand \cite{H}, Jiang \cite{J}, Zhou, Wu and Bian \cite{ZWB}, Zhou and Liu \cite{ZL}, Zhou, Liu and Xu \cite{ZLX}, Sun and Zhou
\cite{SZ}, Zhou, Bian and Pan \cite{ZBP}, Zhou, Wu and Xu \cite{ZWX}, Gao, Guirao and Wu \cite{GGW}, Gao, Guirao and Chen \cite{GGC},
Gao, Wang and Dimitrov \cite{GWD}, Bauer, Nevo and Schmeichel \cite{BNS}, Zhou, Wu and Liu \cite{ZWL}. Zhou \cite{Zb} discussed the
relationship between binding numbers and fractional ID-$k$-factor-critical graphs, and demonstrated a result on a fractional
ID-$k$-factor-critical graph by using a binding number condition of a graph.

\begin{theorem}(\cite{Zb})\label{trivial}
Let $k\geq2$ be an integer, and $G$ be a graph of order $n$ with $n\geq6k-9$. Then $G$ is
fractional ID-$k$-factor-critical if $bind(G)>\frac{(3k-1)(n-1)}{kn-2k+2}$.
\end{theorem}

A graph $G$ is called a fractional $[a,b]$-covered graph if $G$ admits a fractional $[a,b]$-factor with indicator function $h$
satisfying $h(e)=1$ for every $e\in E(G)$. Combining this with the concept of a fractional ID-$[a,b]$-factor-critical graph, it
is natural that we first define the concept of a fractional ID-$[a,b]$-factor-critical covered graph, that is, a graph $G$ is
said to be fractional ID-$[a,b]$-factor-critical covered if $G-I$ is fractional $[a,b]$-covered for any independent set $I$ of
$G$. A fractional ID-$[k,k]$-factor-critical covered graph is simply called a fractional ID-$k$-factor-critical covered graph.
In the previous part of this paper, we introduce the application of the fractional ID-$[a,b]$-factor-critical covered graph.
Now, we recall that the problem on fractional ID-$[a,b]$-factor-critical covered graphs implies that the data packets within
a given capacity range can be still transmitted when certain sites are damaged or blocked, and a channel is assigned in a
communication network, where every site is expressed as a vertex and every channel is modelled as an edge.

Next, we claim a binding number condition for a graph to be fractional ID-$[a,b]$-factor-critical covered, which is a
generalization of Theorem 1.1.

\begin{theorem}\label{trivial}
Let $a$ and $b$ be two integers with $2\leq a\leq b$, and $G$ be a graph of order $n$ with
$n\geq\frac{(a+2b)(a+b-2)+2}{b}$. Then $G$ is fractional ID-$[a,b]$-factor-critical covered if $bind(G)>\frac{(a+2b-1)(n-1)}{bn-(a+b)}$.
\end{theorem}

Naturally, we gain the following result when $a=b=k$ in Theorem 1.2.

\begin{corollary}\label{trivial} Let $k\geq2$ be an integer, and $G$ be a graph of order $n$ with $n\geq6k-4$. Then $G$ is
fractional ID-$k$-factor-critical covered if $bind(G)>\frac{(3k-1)(n-1)}{kn-2k}$.
\end{corollary}

\section{The proof of Theorem 1.2}

Li, Yan and Zhang \cite{LYZ} posed a criterion for a graph being fractional $[a,b]$-covered, which plays a key role in the proof of
Theorem 1.2.

\begin{theorem}(\cite{LYZ})\label{trivial}
Let $a$ and $b$ be two nonnegative integers with $b\geq a$, and $G$ be a graph. Then $G$
is fractional $[a,b]$-covered if and only if for any subset $S\subseteq V(G)$,
$$
\delta_G(S,T)=b|S|-a|T|+d_{G-S}(T)\geq\varepsilon(S),
$$
where $T=\{x:x\in V(G)\setminus S, d_{G-S}(x)\leq a\}$ and $\varepsilon(S)$ is defined by
\[
 \varepsilon(S)=\left\{
\begin{array}{ll}
2,&if \ S \ is \ not \ independent,\\
1,&if \ S \ is \ independent, \ and \ there \ exists \ e=uv\in E(G) \ with \ u\in S,\\
&v\in T \ and \ d_{G-S}(v)=a, \ or \ e_G(S,V(G)\setminus(S\cup T))\geq1,\\
0,&otherwise.\\
\end{array}
\right.
\]
\end{theorem}

Woodall \cite{W} verified the following result, which is also used in the proof of Theorem 1.2.

\begin{lemma}(\cite{W})\label{trivial}
Let $c$ be a positive real number, and let $G$ be a graph of order $n$. If $bind(G)>c$, then
$\delta(G)\geq n-\frac{n-1}{bind(G)}>n-\frac{n-1}{c}$.
\end{lemma}

In what follows, we verify Theorem 1.2.

\begin{proof} Suppose that $G$ satisfies the assumption of Theorem 1.2, but it is not fractional ID-$[a,b]$-factor-critical
covered. Then by Theorem 2.1 and the concept of the fractional ID-$[a,b]$-factor-critical covered graph, there exists some subset $S\subseteq V(H)$
such that
\begin{equation} \label{EQ:1}
\delta_H(S,T)=b|S|-a|T|+d_{H-S}(T)\leq\varepsilon(S)-1,
\end{equation}
where $T=\{x:x\in V(H)\setminus S, d_{H-S}(x)\leq a\}$, $d_{H-S}(T)=\sum\limits_{x\in T}{d_{H-S}(x)}$, $H=G-X$ and $X$ is an independent set of
$G$. In addition, we use $\beta: =bind(G)$ to simplify the notation below.

\medskip
Using Lemma 2.2 and the condition of Theorem 1.2, we gain
\begin{equation} \label{EQ:2}
\delta(G)\geq n-\frac{n-1}{\beta}>\frac{(a+b-1)n+a+b}{a+2b-1}.
\end{equation}

Note that $\varepsilon(S)\leq|S|$. If $T=\emptyset$, then by (\ref{EQ:1}) we possess $\varepsilon(S)-1\geq\delta_H(S,T)=b|S|\geq|S|\geq\varepsilon(S)$, a
contradiction. Hence, $T\neq\emptyset$. Define
$$
h=\min\{d_{H-S}(x):x\in T\}.
$$
From the definition of $T$, we derive that $0\leq h\leq a$.

\medskip
By considering a vertex of $T$, we note that it can possess neighbors in $S$, $X$ and at most $h$ additional neighbors. This gives the following
bound on $\delta(G)$, namely, $\delta(G)\leq|S|+|X|+h$. As a consequence,
\begin{equation} \label{EQ:3}
|S|\geq\delta(G)-|X|-h.
\end{equation}

We now discuss the following two cases.

\medskip
\noindent{\it Case 1.} \ $1\leq h\leq a$.

\smallskip
\noindent Using (\ref{EQ:1}), (\ref{EQ:3}), $|X|\leq n-\delta(G)$ (as $X$ is an independent set), $n\geq|S|+|T|+|X|$ and $\varepsilon(S)\leq2$, we have
\begin{align*}
1\geq&\varepsilon(S)-1\geq\delta_H(S,T)=b|S|-a|T|+d_{H-S}(T)\\
\geq&b|S|-a|T|+h|T|=b|S|-(a-h)|T|\\
\geq&b|S|-(a-h)(n-|X|-|S|)\\
=&(a+b-h)|S|+(a-h)|X|-(a-h)n\\
\geq&(a+b-h)(\delta(G)-|X|-h)+(a-h)|X|-(a-h)n\\
=&(a+b-h)\delta(G)-b|X|-h(a+b-h)-(a-h)n\\
\geq&(a+b-h)\delta(G)-b(n-\delta(G))-h(a+b-h)-(a-h)n\\
=&(a+2b-h)\delta(G)-(a+b-h)n-h(a+b-h).
\end{align*}
Solving for $\delta(G)$, we derive the following
\begin{align*}
\delta(G)\leq f(h): =\frac{(a+b-h)(n+h)+1}{a+2b-h}.
\end{align*}

Taking the derivative of $f(h)$ with respect to $h$ yields
\begin{align*}
\frac{df}{dh}=&\frac{(a+2b-h)(-(n+h)+(a+b-h))+((a+b-h)(n+h)+1)}{(a+2b-h)^{2}}\\
=&\frac{-bn+a^{2}+3ab+2b^{2}-2ah-4bh+h^{2}+1}{(a+2b-h)^{2}}\\
\leq&\frac{-bn+a^{2}+3ab+2b^{2}-2a-4b+1+1}{(a+2b-h)^{2}}\\
=&\frac{-bn+(a+2b)(a+b-2)+2}{(a+2b-h)^{2}}.
\end{align*}

For $n\geq\frac{(a+2b)(a+b-2)+2}{b}$, we derive that $\frac{df}{dh}\leq0$, implying that $f(h)$ attains its maximum at smallest value of $h$.
Therefore,
\begin{align*}
\delta(G)\leq\frac{(a+b-1)(n+1)+1}{a+2b-1}=\frac{(a+b-1)n+a+b}{a+2b-1},
\end{align*}
this contradicts (\ref{EQ:2}).

\medskip
\noindent{\it Case 2.} \ $h=0$.
\smallskip\\
\noindent{\it Subcase 2.1.} \ $\beta\leq a+b-1$.
\smallskip\\
Setting $Z=\{x:x\in T,d_{H-S}(x)=0\}$. Evidently, $Z\neq\emptyset$ and $N_G(V(G)\setminus(X\cup S))\cap Z=\emptyset$, which hints
$|N_G(V(G)\setminus(X\cup S))|\leq n-|Z|$. Thus,
\begin{align*}
bind(G)=\beta\leq\frac{|N_G(V(G)\setminus(X\cup S))|}{|V(G)\setminus(X\cup S)|}\leq\frac{n-|Z|}{n-|X|-|S|},
\end{align*}
namely,
\begin{equation} \label{EQ:4}
|S|\geq\left(1-\frac{1}{\beta}\right)n-|X|+\frac{1}{\beta}|Z|.
\end{equation}

Using (\ref{EQ:1}), (\ref{EQ:2}), (\ref{EQ:4}), $2\leq a\leq b$, $Z\neq\emptyset$, $|X|\leq n-\delta(G)$, $n\geq|S|+|T|+|X|$ and $\varepsilon(S)\leq2$,
we acquire
\begin{align*}
1\geq&\varepsilon(S)-1\geq\delta_H(S,T)=b|S|-a|T|+d_{H-S}(T)\\
\geq&b|S|-a|T|+|T|-|Z|\\
=&b|S|-(a-1)|T|-|Z|\\
\geq&b|S|-(a-1)(n-|X|-|S|)-|Z|\\
=&(a+b-1)|S|-(a-1)n+(a-1)|X|-|Z|\\
\geq&(a+b-1)\left(\left(1-\frac{1}{\beta}\right)n-|X|+\frac{1}{\beta}|Z|\right)-(a-1)n+(a-1)|X|-|Z|\\
=&bn-\frac{a+b-1}{\beta}n-b|X|+\left(\frac{a+b-1}{\beta}-1\right)|Z|\\
\geq&bn-\frac{a+b-1}{\beta}n-b|X|+\left(\frac{a+b-1}{\beta}-1\right)\\
\geq&bn-\frac{a+b-1}{\beta}n-b(n-\delta(G))+\frac{a+b-1}{\beta}-1\\
=&-\frac{a+b-1}{\beta}n+b\delta(G)+\frac{a+b-1}{\beta}-1\\
\geq&-\frac{a+b-1}{\beta}n+b\left(n-\frac{n-1}{\beta}\right)+\frac{a+b-1}{\beta}-1\\
\geq&-\frac{a+b-1}{\beta}n+b\left(n-\frac{n-1}{\beta}\right)+\frac{a+b-1}{\beta}-(a+b-1)\\
=&-\frac{(a+2b-1)n-b}{\beta}+bn-(a+b-1)\left(1-\frac{1}{\beta}\right).
\end{align*}
Solving for $\beta$, this yields:
\begin{align*}
\beta\leq\frac{(a+2b-1)(n-1)}{bn-(a+b)},
\end{align*}
which contradicts the condition of Theorem 1.2.

\medskip
\noindent{\it Subcase 2.2.} \ $\beta>a+b-1$.
\smallskip\\
Applying (\ref{EQ:2}) and $2\leq a\leq b$, we achieve
\begin{equation} \label{EQ:5}
\delta(G)\geq n-\frac{n-1}{\beta}>n-\frac{n-1}{a+b-1}=\frac{(a+b-2)n+1}{a+b-1}\geq\frac{(a+b)n}{2a+b}+\frac{1}{a+b-1}.
\end{equation}

It follows from (\ref{EQ:1}), (\ref{EQ:3}), (\ref{EQ:5}), $2\leq a\leq b$, $|X|\leq n-\delta(G)$, $n\geq|S|+|T|+|X|$ and $\varepsilon(S)\leq2$ that
\begin{align*}
1\geq&\varepsilon(S)-1\geq\delta_H(S,T)=b|S|-a|T|+d_{H-S}(T)\\
\geq&b|S|-a|T|\geq b|S|-a(n-|X|-|S|)\\
=&(a+b)|S|-an+a|X|\\
\geq&(a+b)(\delta(G)-|X|)-an+a|X|\\
=&(a+b)\delta(G)-an-b|X|\\
\geq&(a+b)\delta(G)-an-b(n-\delta(G))\\
=&(a+2b)\delta(G)-(a+b)n\\
>&(a+2b)\left(\frac{(a+b)n}{2a+b}+\frac{1}{a+b-1}\right)-(a+b)n\\
=&\frac{a+2b}{a+b-1}>1,
\end{align*}
a contradiction. We certify Theorem 1.2.
\end{proof}

\section{Conclusion}

In this work, we demonstrate a binding number condition for a graph to be fractional ID-$[a,b]$-factor-critical covered. But, we do
not know whether the bound on $bind(G)$ in Theorem 1.2 is sharp or not. Naturally, we put forward the following conjecture:

\medskip

\noindent{\textbf{Conjecture 3.1.}} Let $a$ and $b$ be two integers with $2\leq a\leq b$, and $G$ be a graph of order $n$ with
$n\geq\frac{(a+2b)(a+b-2)+2}{b}$. Then $G$ is fractional ID-$[a,b]$-factor-critical covered if $bind(G)\geq\frac{(a+2b-1)(n-1)}{bn-(a+b)}$.

\medskip

In the proof of Theorem 1.2, the condition $bind(G)>\frac{(a+2b-1)(n-1)}{bn-(a+b)}$ is necessary. But for Conjecture 3.1, I do not know how to
prove it. Next, we argue the extreme case of $a=b=k$, then a fractional ID-$[a,b]$-factor-critical covered graph is a fractional ID-$k$-factor-critical
covered graph, which is an extension of a fractional ID-$k$-factor-critical graph. And so, Theorem 1.2 in this paper is a generalization of
Zhou's previous result \cite{Zb}. Furthermore, we introduce the applications of the fractional $[a,b]$-factors of graphs and the fractional
ID-$[a,b]$-factor-critical covered graphs in Section~1.

\subsection*{Acknowledgments}

I take this opportunity to thank the anonymous referees for their careful reading of the manuscript and suggestions
which have immensely helped us in getting the paper to its present form.

\section*{Declaration of interest statement}

The author declares that there is no conflict of interests regarding the publication of this paper.


\begin{thebibliography}{10}
\providecommand{\url}[1]{\texttt{#1}}
\providecommand{\urlprefix}{URL }
\expandafter\ifx\csname urlstyle\endcsname\relax
  \providecommand{\doi}[1]{doi:\discretionary{}{}{}#1}\else
  \providecommand{\doi}{doi:\discretionary{}{}{}\begingroup
  \urlstyle{rm}\Url}\fi
\providecommand{\eprint}[2][]{\url{#2}}

\bibitem{AMB}
Araujo DRB, Martins JF, Bastos CJA.
\newblock New graph model to design optical networks.
\newblock \emph{IEEE Communications Letters}, 2015.
\newblock \textbf{19}(12):2130--2133.
\newblock \doi{10.1109/LCOMM.2015.2480716}.

\bibitem{AP}
Ashwin P, Postlethwaite C.
\newblock On designing heteroclinic networks from graphs.
\newblock \emph{Physica D}, 2013.
\newblock \textbf{265}:26--39.
\newblock \doi{10.1016/j.physd.2013.09.006}.

\bibitem{BNS}
Bauer D, Nevo A, Schmeichel E.
\newblock Best monotone degree condition for the Hamiltonicity of graphs with a 2-factor.
\newblock \emph{Graphs and Combinatorics}, 2017.
\newblock \textbf{33}(5):1231--1248.
\newblock \doi{10.1007/s00373-017-1840-1}.

\bibitem{FLJ}
Fardad M, Lin F, Jovanovic MR.
\newblock Design of optimal sparse interconnection graphs for synchronization of oscillator networks.
\newblock \emph{IEEE Transactions on Automatic Control}, 2014.
\newblock \textbf{59}(9):2457--2462.
\newblock \doi{10.1109/TAC.2014.2301577}.

\bibitem{GGW}
Gao W, Guirao J, Wu H.
\newblock Two tight independent set conditions for fractional $(g,f,m)$-deleted graphs systems.
\newblock \emph{Qualitative Theory of Dynamical Systems}, 2018.
\newblock \textbf{17}(1):231--243.
\newblock \doi{10.1007/s12346-016-0222-z}.

\bibitem{GGC}
Gao W, Guirao J, Chen Y.
\newblock A toughness condition for fractional $(k,m)$-deleted graphs revisited.
\newblock \emph{Acta Mathematica Sinica, English Series}, 2019.
\newblock \textbf{35}(7):1227--1237.
\newblock \doi{10.1007/s10114-019-8169-z}.

\bibitem{GWD}
Gao W, Wang W, Dimitrov D.
\newblock Toughness condition for a graph to be all fractional $(g,f,n)$-critical deleted.
\newblock \emph{Filomat}, 2019.
\newblock \textbf{33}(9):2735--2746.
\newblock \doi{10.2298/FIL1909735G}.

\bibitem{HK}
Haghparast N, Kiani D.
\newblock Edge-connectivity and edges of even factors of graphs.
\newblock \emph{Discussiones Mathematicae Graph Theory}, 2019.
\newblock \textbf{39}(2):357--364.
\newblock \doi{10.7151/dmgt.2082}.

\bibitem{H}
Hasanvand M.
\newblock Factors and connected factors in tough graphs with high isolated toughness.
\newblock \emph{Dec. 30. arXiv: 1812.11640}, 2018.
\newblock \doi{10.48550/arXiv.1812.11640}.

\bibitem{J}
Jiang J.
\newblock A sufficient condition for all fractional $[a,b]$-factors in graphs.
\newblock \emph{Proceedings of the Romanian Academy, Series A}, 2018.
\newblock \textbf{19}(2):315--319.

\bibitem{LMA}
Lanzeni S, Messina E, Archetti F.
\newblock Graph models and mathematical programming in biochemical network analysis and metabolic engineering design.
\newblock \emph{Computers and Mathematics with Applications}, 2008.
\newblock \textbf{55}(5):970--983.
\newblock \doi{10.1016/j.camwa.2006.12.101}.

\bibitem{LYZ}
Li Z, Yan G, Zhang X.
\newblock On fractional $(g,f)$-covered graphs.
\newblock \emph{OR Transactions (China)}, 2002.
\newblock \textbf{6}(4):65--68.

\bibitem{PR}
Pishvaee MS, Rabbani M.
\newblock A graph theoretic--based heuristic algorithm for responsive supply chain network design with direct and indirect shipment.
\newblock \emph{Advances in Engineering Software}, 2011.
\newblock \textbf{42}(3):57--63.
\newblock \doi{10.1016/j.advengsoft.2010.11.001}.

\bibitem{RH}
Rahimi M, Haghighi A.
\newblock A graph portioning approach for hydraulic analysis-design of looped pipe networks.
\newblock \emph{Water Resources Management}, 2015.
\newblock \textbf{29}(14):5339--5352.
\newblock \doi{10.1007/s11269-015-1121-9}.

\bibitem{SZ}
Sun Z, Zhou S.
\newblock A generalization of orthogonal factorizations in digraphs.
\newblock \emph{Information Processing Letters}, 2018.
\newblock \textbf{132}:49--54.
\newblock \doi{10.1016/j.ipl.2017.12.003}.

\bibitem{W}
Woodall DR.
\newblock The binding number of a graph and its Anderson number.
\newblock \emph{Journal of Combinatorial Theory, Series B}, 1973.
\newblock \textbf{15}(3):225--255.
\newblock \doi{10.1016/0095-8956(73)90038-5}.

\bibitem{WZi}
Wang S, Zhang W.
\newblock Isolated toughness for path factors in networks.
\newblock \emph{RAIRO-Operations Research}, 2022.
\newblock \textbf{56}(4):2613--2619.
\newblock \doi{10.1051/ro/2022123}.

\bibitem{WZo}
Wang S, Zhang W.
\newblock On $k$-orthogonal factorizations in networks.
\newblock \emph{RAIRO-Operations Research}, 2021.
\newblock \textbf{55}(2):969--977.
\newblock \doi{10.1051/ro/2021037}.

\bibitem{WZr}
Wang S, Zhang W.
\newblock Research on fractional critical covered graphs.
\newblock \emph{Problems of Information Transmission}, 2020.
\newblock \textbf{56}(3):270--277.
\newblock \doi{10.1134/S0032946020030047}.

\bibitem{Za}
Zhou S.
\newblock A neighborhood union condition for fractional $(a,b,k)$-critical covered graphs.
\newblock \emph{Discrete Applied Mathematics}, 2021.
\newblock \doi{10.1016/j.dam.2021.05.022}.

\bibitem{Za0}
Zhou S.
\newblock A result on fractional $(a,b,k)$-critical covered graphs.
\newblock \emph{Acta Mathematicae Applicatae Sinica-English Series}, 2021.
\newblock \textbf{37}(4):657--664.
\newblock \doi{10.1007/s10255-021-1034-8}.

\bibitem{Zb}
Zhou S.
\newblock Binding numbers for fractional ID-$k$-factor-critical graphs.
\newblock \emph{Acta Mathematica Sinica, English Series}, 2014.
\newblock \textbf{30}(1):181--186.
\newblock \doi{10.1007/s10114-013-1396-9}.

\bibitem{Zr}
Zhou S.
\newblock Remarks on restricted fractional $(g,f)$-factors in graphs.
\newblock \emph{Discrete Applied Mathematics}, 2022.
\newblock \doi{10.1016/j.dam.2022.07.020}.

\bibitem{ZB}
Zhou S, Bian Q.
\newblock The existence of path-factor uniform graphs with large connectivity.
\newblock \emph{RAIRO-Operations Research}, 2022.
\newblock \doi{10.1051/ro/2022143}.

\bibitem{ZBP}
Zhou S, Bian Q, Pan Q.
\newblock Path factors in subgraphs.
\newblock \emph{Discrete Applied Mathematics}, 2022.
\newblock \textbf{319}:183--191.
\newblock \doi{10.1016/j.dam.2021.04.012}.

\bibitem{ZL}
Zhou S, Liu H.
\newblock Discussions on orthogonal factorizations in digraphs.
\newblock \emph{Acta Mathematicae Applicatae Sinica-English Series}, 2022.
\newblock \textbf{38}(2):417--425.
\newblock \doi{10.1007/s10255-022-1086-4}.

\bibitem{ZLX}
Zhou S, Liu H, Xu Y.
\newblock A note on fractional ID-$[a,b]$-factor-critical covered graphs.
\newblock \emph{Discrete Applied Mathematics}, 2022.
\newblock \textbf{319}:511--516.
\newblock \doi{10.1016/j.dam.2021.03.004}.

\bibitem{ZSB}
Zhou S, Sun Z, Bian Q.
\newblock Isolated toughness and path-factor uniform graphs (II).
\newblock \emph{Indian Journal of Pure and Applied Mathematics}, 2022.
\newblock \doi{10.1007/s13226-022-00286-x}.

\bibitem{ZWB}
Zhou S, Wu J, Bian Q.
\newblock On path-factor critical deleted (or covered) graphs.
\newblock \emph{Aequationes Mathematicae}, 2022.
\newblock \textbf{96}(4):795--802.
\newblock \doi{10.1007/s00010-021-00852-4}.

\bibitem{ZWL}
Zhou S, Wu J, Liu H.
\newblock Independence number and connectivity for fractional $(a,b,k)$-critical covered graphs.
\newblock \emph{RAIRO-Operations Research}, 2022.
\newblock \textbf{56}(4):2535--2542.
\newblock \doi{10.1051/ro/2022119}.

\bibitem{ZWX}
Zhou S, Wu J, Xu Y.
\newblock Toughness, isolated toughness and path factors in graphs.
\newblock \emph{Bulletin of the Australian Mathematical Society}, 2021.
\newblock \doi{10.1017/S0004972721000952}.

\end{thebibliography}


\end{document}